\documentclass[12pt]{article}

\usepackage[a4paper,text={17.5cm,24.0cm},centering]{geometry}
\usepackage{latexsym,graphicx,amsmath,amssymb,mathrsfs,bm}
\usepackage{fancybox,color,wrapfig,framed,mathtools}
\usepackage[UKenglish]{babel}
\def\bM{{\bf M}}
\def\bK{{\bf K}}
\def\bJ{{\bf J}}
\def\ri{\texttt{i}}
\def\rd{{\rm d}}
\def\re{{\rm e}}
\def\fr{\mbox{$\frac{1}{2}$}}

\def\R{\mathbb R}

\def\qand{\quad\mbox{and}\quad}
\def\bq{{\bf q}}
\def\bp{{\bf p}}
\def\br{{\bf r}}
\def\bu{{\bf u}}
\def\bw{{\bf w}}
\def\bk{{\bf k}}
\def\eps{{\varepsilon}}
\def\hatu{\widehat{\bf u}}

\def\cA{\mathscr{A}}
\def\cB{\mathscr{B}}
\def\Zh{\widehat Z}
\def\Rg{{\mathcal R}}

\def\langang{\langle\!\langle}
\def\rangang{\rangle\!\rangle}

\begin{document}

\begin{center}
  \textbf{\Large Regularizing the Cross-Newell equation}\\[2mm]
  \textbf{\Large by re-modulating along a characteristic angle}
\vspace{.75cm}

\textsf{\large Nicholas J. Burgess and Thomas J. Bridges}
\vspace{.25cm}

\textit{Department of Mathematics, University of Surrey, Guildford GU2 7XH, UK}\\
Email: \textcolor{blue}{\texttt{N.Burgess@surrey.ac.uk}} {\small and}
\textcolor{blue}{\texttt{T.Bridges@surrey.ac.uk}}
\vspace{1.0cm}

\setlength{\fboxsep}{10pt}
\doublebox{\parbox{15cm}{
    {\bf Abstract.} {\small
      A regularization of the Cross-Newell equation is
      presented.  It is based on a secondary
      re-modulation along characteristics. This new characteristic Cross-Newell
      equation is not isotropic (has preferred directions), but
      is universal (homogeneous in $\eps$
      and equation independent coefficients), has fourth derivatives,
      and generates localized
      solutions that are bi-asymptotic to rolls. Re-modulation of rolls in the
      Ginzburg-Landau equation is used for illustration of the theory.}
}}
\end{center}

\section{Introduction}
\setcounter{equation}{0}
\label{sec-intro}

The Cross-Newell equation is a second-order partial differential
equation that
captures the modulation, phase singularities, bifurcations, and
stability of
steady periodic patterns in the plane.  It is derived via
reduction, in the neighborhood of the family of periodic patterns,
from a class of nonlinear partial differential equations
(e.g.\ \textsc{Cross \& Newell}~\cite{cn84}, \textsc{Ercolani, et al.}~\cite{einp00}, \textsc{Hoyle}~\cite{hoyle},
\textsc{Cross \& Greenside}~\cite{cg09}, and references therein).
In standard form it is
\begin{equation}\label{cn-eqn}
\tau(q)\phi_T = -\nabla\cdot\big(\mathcal{B}(q){\bf q}\big) \,,
\end{equation}
where $\bq=(q_1,q_2)=\nabla\phi$
is the vector-valued modulation of the wavenumber and $q=\|\bq\|$.
The functions $\tau(q)>0$ and $\mathcal{B}$ are determined from
properties of the governing equations, and they have a universal form.
The independent variables here are scaled with
\begin{equation}\label{cn-scale}
  X=\eps x\,,\quad Y=\eps y\,,\quad T=\eps^2 t\,,
\end{equation}
where $0<\eps\ll 1$ is a small
parameter with $\eps^{-1}$ a measure of the size of the domain.

The \emph{regularized Cross-Newell equation} is obtained by taking the
perturbation expansion to the next order
\begin{equation}\label{rcn-eqn}
  \tau(q)\phi_T = -\nabla\cdot\big(\mathcal{B}(q){\bf q}\big) -\eps^2\eta
  \Delta^2\phi\,,
\end{equation}
(e.g.\ \S2.5 of \cite{einp00}).  Here $\eta>0$ is deduced from the
governing equations and $\Delta$ is the Laplacian in $(X,Y)$.
A key advantage of this regularization, in addition to smoothing, is
that it is isotropic.  A disadvantage is that it is not homogeneous in
$\eps$ and so reduces to the Cross-Newell equation in the limit $\eps\to0$.

The Cross-Newell (CN) equation is invariant under rotations in the plane in
the following sense.  Introducing rotated coordinates
\[
Z_1 = X+c_1 Y \qand Z_2 = X+c_2Y\,,\quad c_1\neq c_2\,,
\]
in (\ref{cn-eqn}) leaves the form of the equation invariant, changing
only the formulae for the coefficients.  There is no preferred direction
in the CN equation.
On the other hand there
are special directions in the plane, the characteristic directions,
where the CN equation can be re-modulated to generate fourth derivative
terms.  

The purpose of this paper is to follow these directions and
introduce a new ``regularization'' to the
CN equation. This ``characteristic Cross-Newell'' (CCN) equation
has a preferred direction and so is not isotropic, but is universal
(regularizing term is independent of $\eps$, equation-independent
coefficients), generates localized solutions,
and is generic (exists for an open set of wavenumbers).
The principal disadvantages are twofold: ``regularization'' is in
quotes as the higher order term may be destabilizing, and secondly
the equation is valid in the zig-zag unstable regime only.

The CCN equation is
\begin{equation}\label{nrcn-eqn}
  \tau(q)\phi_T = \pm 2\sqrt{-\Delta_{zz}} \,\phi_{XY}
  + \kappa \phi_{X}\phi_{XX} + \mathscr{K} \phi_{XXXX} \,,
\end{equation}
where $\Delta_{zz}$ is defined below, and
$\phi$ is related to the modulation wavenumber as
$\bq = \nabla\phi$ with $\nabla$ in $(X,Y)$ variables.  The
independent variables are re-scaled to assure that
(\ref{nrcn-eqn}) is homogeneous in $\eps$, and
are skewed to align with a characteristic direction
\begin{equation}\label{XYT-def}
X = \eps(x + C^{\pm} y)\,,\quad Y=\eps^3 y\,,\quad\mbox{and}\quad
T = \eps^4t\,.
\end{equation}
Here $C^{\pm}$ are characteristics of the steady CN equation.  There
are two CCN equations, one associated with each of the characteristics
$C^\pm$.
The $\phi_{XY}$ term in (\ref{nrcn-eqn}) is a diagonalization of the
linearized CN equation (\ref{steady-cn-linear}) and this is shown in
\S\ref{sec-char-cn}.
The coefficient of ``regularization'', $\mathscr{K}$, may be
positive or negative. The significance and origin
of the coefficients in (\ref{nrcn-eqn}) will be discussed in
\S\ref{rccn-eqn-derivation}.
The most novel is $\kappa$ which can be expressed in terms of
the fluxes of a conservation law with components $(\cB,\cA)$,
\begin{equation}\label{kappa-def}
\kappa = \left(\frac{\partial\ }{\partial k}+ C^{\pm}\frac{\partial\ }{\partial\ell} \right)^2\big(
\cB + C^{\pm}\cA\big)\,,
\end{equation}
where ${\bf k}=(k,\ell)$ is the wavenumber of the basic state.
The conservation
law is due to a symmetry and is introduced in \S\ref{subsec-claw}.  

The equation (\ref{nrcn-eqn}) is obtained by re-modulation once the
characteristics of the original CN equation are known.
The characteristic directions $C^{\pm}$ are obtained from the
linearized CN equation,
also called the phase diffusion equation, which can be written in the form
\begin{equation}\label{steady-cn-linear}
  \tau\phi_T =
  \cB_k \phi_{XX} + (\cB_\ell+\cA_k)\phi_{XY} +
  \cA_{\ell}\phi_{YY}  \,,
\end{equation}
where here we have reverted to the CN scaling (\ref{cn-scale}).
The coefficients on the right-hand side are based
on the $(\cB,\cA)$ conservation law (as in $\kappa$ in (\ref{kappa-def}),
and discussed in \S\ref{subsec-claw}).
For now they are just $(k,\ell)-$dependent scalars.  Let
\[
\phi = \widehat\phi \re^{\ri(X+CY)}
\]
then the characteristics of the steady CN equation,
denoted by $C$, satisfy the quadratic
\[
\cA_{\ell} C^2 + (\cB_\ell+\cA_k) C +
  \cB_k = 0\,,
  \]
  with roots
  \begin{equation}\label{C-plusminus-def}
  C^{\pm} = -\frac{(\cB_\ell+\cA_k)}{2\cA_{\ell}}
    \pm \frac{1}{\cA_{\ell}}\sqrt{-\Delta_{zz}}\,,
  \end{equation}
    where
    \begin{equation}\label{Delta-zz-def}
    \Delta_{zz}:=
    {\rm det}\left[\begin{matrix} 
          \cA_{\ell} & \cA_{k} \\
          \cB_{\ell} & \cB_{k}\end{matrix}\right]\,.
    \end{equation}
    It is shown in \S\ref{sec-char-cn} below that
    a negative (positive) sign of $\Delta_{zz}$ signals a zig-zag
    instability (stability). The above formula shows that
    zig-zag instability is correlated with real characteristics.
    The region of existence of rolls, denoted by $\mathcal{D}(k,\ell)\subset\R^2$,
    splits into two regions $\mathcal{D}=\mathcal{D}^-\cup\mathcal{D}^+$,
    with $\mathcal{D}^{\pm}$ the subset of $\mathcal{D}$ where
    ${\rm sign}(\Delta_{zz})=\pm1$.

    Henceforth we will assume that the characteristics are real:
    $(k,\ell)\in\mathcal{D}^-$.
    In principle the theory goes through when the basic state is
    zig-zag stable and the characteristics are
    complex; that is, $(k,\ell)\in\mathcal{D}^+$.
    But then the independent variables and coefficients in the
    resulting CCN modulation equation are complex. 
    
    This strategy of using characteristic directions to re-modulate
    is motivated by the
recent discovery of \textsc{Ratliff}~\cite{r19} that Whitham modulation
theory always generates dispersion on a long enough time scale when the
characteristics are hyperbolic.  The emergence of dispersion is shown
there by re-modulating Whitham theory along characteristic directions.

Although the equation (\ref{nrcn-eqn}) is not isotropic, there
are three directions associated with the theory.
Firstly there is the angle $\vartheta$ induced by the choice of wavenumber
vector
\begin{equation}\label{vartheta-def}
  {\bf k}=(k,\ell)= \|\bk\|(\cos\vartheta,\sin\vartheta)\,,
\end{equation}
  of the basic state.
Secondly there are the characteristic
angles induced by $C^{\pm}$ along which the
equation (\ref{nrcn-eqn}) is operational.
Thirdly, the localized solutions of the steady
version of (\ref{nrcn-eqn}) will be aligned with a third angle
(cf.\ \S\ref{sec-solitons}).

Our strategy for deriving (\ref{nrcn-eqn})
is similar to the CN derivation, but with
a new scaling valid along characteristic directions.
The starting point is an abstract set of nonlinear PDEs,
\begin{equation}\label{LuNu}
  {\bf u}_t = {\bf Lu} + \mathcal{N}(\bu) \,,
\end{equation}
for some vector-valued unknown $\bu(x,y,t)$, linear operator ${\bf L}$
and nonlinear operator $\mathcal{N}$.
However, in order
to achieve universality, whereby the coefficients in
(\ref{nrcn-eqn}) are deduced from abstract properties of the
governing equations rather than each specific equation,
additional structure of the governing equations
is necessary.  Firstly, the right-hand side
of (\ref{LuNu}) is required to be of gradient type.
Secondly, we use an exact conservation law for steady periodic patterns
in the plane.  It is
the fluxes of this conservation law that generate $\cA$ and
$\cB$ in (\ref{kappa-def}) and (\ref{steady-cn-linear}).

Assuming a basic
periodic pattern exists,
\begin{equation}\label{basicstate-1}
\bu(x,y,t) = \hatu(\theta,{\bf k})\,,\quad \theta = k x+\ell y + \theta_0\,,
\end{equation}
$2\pi-$periodic in $\theta$,
the ansatz used in the Cross-Newell derivation is,
\begin{equation}\label{ansatz-1}
\bu(x,y,t) = \hatu(\theta+ \eps^{-1}\phi,
   {\bf k}+\bq) + \eps {\bf w}(\theta+\eps^{-1}\phi,X,Y,T,\eps)\,,
\end{equation}
with $X=\eps x$, $Y=\eps y$, and $T=\eps^2t$, and $\bq=\nabla\phi$ in
these variables, is replaced by
\begin{equation}\label{ansatz-2}
\bu(x,y,t) = \hatu(\theta+ \eps\phi,
   {\bf k}+\eps^2{\bf C}(\eps){\bf q}) + \eps^3 {\bf w}(\theta+\eps\phi,X,Y,T,\eps)\,,
\end{equation}
but now with the characteristic scaled variables (\ref{XYT-def}),
and $\bq=\nabla\phi$ also in the variables (\ref{XYT-def}).
The precise form of the
matrix ${\bf C}(\eps)$ will be introduced in \S\ref{rccn-eqn-derivation}.
   In both (\ref{ansatz-1}) and (\ref{ansatz-2})
   the term ${\bf w}$ is a remainder term. Expansion and substitution
  of these ans\"atze into the governing equations, and solvability, leads to
  the CN equation (\ref{cn-eqn}) in the case of (\ref{ansatz-1})
  and to the CCN equation (\ref{nrcn-eqn}) in the case of (\ref{ansatz-2}).

  For illustration of the theory, it is applied to the real
  Ginzburg-Landau equation,
  \begin{equation}\label{gl-real}
  \Psi_t = \Psi_{xx}+\Psi_{yy} + \Psi - |\Psi|^2 \Psi\,,
  \end{equation}
  for complex-valued $\Psi(x,y,t)$. Although this example is elementary,
  it contains all the key features of the theory.  The basic family of
  rolls, and the functions $\cB(k,\ell)$ and $\cA(k,\ell)$ can be
  obtained explicitly, and the coefficients $\kappa$ and
  $\mathscr{K}$ can be computed explicitly.  It is found that $\mathscr{K}<0$
  for all admissable $(k,\ell)$, and so the CCN equation provides
  a proper regularization of
  the CN equation along characteristics.
  
  An outline of the paper is as follows.  First the linear equation
  (\ref{steady-cn-linear}) is analyzed in \S\ref{sec-char-cn} for
    characteristics, zig-zag instability, and the implications of re-scaling.
    Then in \S\ref{sec-ge} the class of
    governing equations and planar conservation law
    are introduced.  A brief look at the derivation of the classical
    CN equation is given in \S\ref{sec-cn-eqn-sketch}
    to introduce the role of the conservation
    law, and then in \S\ref{rccn-eqn-derivation}
    a derivation of (\ref{nrcn-eqn}) is given.
    In \S\ref{sec-gl} the theory is applied to the real Ginzburg-Landau equation
    (\ref{gl-real}).
    Finally in the concluding remarks section some implications and
    other directions are discussed.

    \section{The characteristic angles}
    \label{sec-char-cn}
    \setcounter{equation}{0}

    The linear CN equation (\ref{steady-cn-linear})
    can be diagonalized two ways
\begin{equation}\label{cn-linear-diag1}
  \tau\phi_T =
  \cB_k\left( \frac{\partial\ }{\partial X} + \gamma_1\frac{\partial\ }{\partial Y}\right)^2\phi + \frac{1}{\cB_k}\Delta_{zz} \frac{\partial^2\phi}{\partial Y^2}\,,\quad \gamma_1=\frac{(\cB_\ell+\cA_k)}{2\cB_k}\,,
\end{equation}
and
\begin{equation}\label{cn-linear-diag2}
  \tau\phi_T =\cA_\ell\left( \frac{\partial\ }{\partial Y} + \gamma_2\frac{\partial\ }{\partial X}\right)^2\phi + \frac{1}{\cA_{\ell}}\Delta_{zz}
  \frac{\partial^2\phi}{\partial X^2}\,,
  \quad \gamma_1=\frac{(\cB_\ell+\cA_k)}{2\cA_\ell}\,.
\end{equation}
A negative sign of either $\cA_\ell$ or $\cB_k$ is associated
with Eckhaus instability.  The Eckhaus instability or stability
does not enter the theory in any way.  On the other hand, regardless
of the sign of $\cA_\ell$ or $\cB_k$, it is clear that
the linear system (\ref{steady-cn-linear})
is ill-posed when $\Delta_{zz}<0$; that is, the zig-zag instability criterion.
Ill-posedness of (\ref{steady-cn-linear}) indicates that the
basic periodic pattern is unstable.
On the other hand, as shown in (\ref{C-plusminus-def}),
$\Delta_{zz}<0$ is also the condition for real characteristics.  Henceforth
it is assumed that the characteristics are real
\begin{equation}\label{zig-zag-unstable}
  \Delta_{zz}<0 \quad\Rightarrow\quad (k,\ell)\in\mathcal{D}^-\,.
\end{equation}
Our strategy is to remodulate the basic state using
the characteristic coordinates (\ref{XYT-def}).  To see the
impact of the characteristic coordinates
on the linear CN equation (\ref{steady-cn-linear}), let
\begin{equation}\label{Z1Z2-scaling}
Z_1 = X+C^{\pm}Y = \eps(x + C^{\pm} y)\qand Z_2=\eps^2 Y=\eps^3y\,,
\end{equation}
and substitute into (\ref{steady-cn-linear}),
\begin{equation}\label{Z1Z2-eqn}
\tau \phi_T = \pm 2 \eps^2 \sqrt{-\Delta_{zz}} \frac{\partial^2\phi}{\partial Z_1\partial Z_2} + \mathcal{O}(\eps^4)\,.
\end{equation}
The new scaling (\ref{Z1Z2-scaling}) drives the pure second derivatives to
higher order, leaving the cross-derivative term only.  When $Z_1$ and
$Z_2$ are replaced by the symbols (\ref{XYT-def}), we see how the first term in
(\ref{nrcn-eqn}) is generated.  For consistency it is clear that
the scaling $T=\eps^2$ should be modified to $T=\eps^4t$.
The equation (\ref{Z1Z2-eqn})
also indicates that the nonlinearity and the fourth order
term in (\ref{nrcn-eqn}) should be of order $\eps^4$.
As we will see, all this will work with the right choice
of ansatz (\ref{ansatz-2}).

    \section{The governing equations}
    \label{sec-ge}
    \setcounter{equation}{0}
    
    The governing equations can be taken in the abstract form (\ref{LuNu}) with
    the requirement that the right hand side is a gradient 
    \[
    \bu_t = \frac{\delta \mathcal{L}}{\delta\bu}\,,
    \]
    for some functional $\mathcal{L}(\bu,\bu_x,\bu_y)$.  However we will
    go one step further and Legendre transform $\mathcal{L}$ in both the
    $\bu_x$ and $\bu_y$ directions.  This leads to a multisymplectic gradient
    structure.  Examples of the transformation of functionals into
    multisymplectic form are given in \cite{bhl10} and \cite{tjb17}.
    We take as a starting point that $\mathcal{L}$ is so transformed,
    with new vector-valued dependent variables $Z(x,y,t)$,
    \begin{equation}\label{eq1}
      \mathcal{L} = \iint_\mathcal{D}\left[
        \fr \langle \textbf{J} Z_x,Z \rangle
      + \fr \langle \textbf{K}Z_y, Z \rangle - S(Z) \right]\,\rd x\rd y\,,
    \end{equation}
    where $Z\in\R^{2n}$, $\mathcal{D}$ is an open subset of $\R^2\ni(x,y)$,
    and $\langle\cdot,\cdot\rangle$ is the inner product on $\R^{2n}$.
    The matrices $\bJ$ and $\bK$ are skew symmetric, and $S\,:\,\R^{2n}\to\R$
    is a generalized Hamiltonian function.  Taking the variational derivative
    and bringing time back in gives the following canonical form for the
    governing equations.
    \begin{equation}\label{MJKS-eqn}
      \bM Z_t + \bJ Z_x + \bK Z_y = \nabla S(Z)\,.
    \end{equation}
    In this equation $\bM$ is symmetric and positive semi-definite.
    In \S\ref{sec-gl} the real Ginzburg-Landau (\ref{gl-real})
    is transformed to the structural form (\ref{MJKS-eqn}).
    
\subsection{Symmetry in the plane and a conservation law}
\label{subsec-claw}

The steady equations have a conservation law which follows from the
gradient structure of the right hand side of (\ref{LuNu}) and
\eqref{MJKS-eqn}.  This conservation law is reminiscent of
conservation of wave action in Whitham theory
(cf.\ \textsc{Newell \& Pomeau}~\cite{np95}).
When the equations are in the form (\ref{MJKS-eqn}) the
conservation law can be given a tidy geometric form.
Consider a one-parameter family of steady solutions,
$Z(x,y,s)$, of (\ref{MJKS-eqn})
parameterized by $s$, an ensemble parameter, and suppose it is a loop
of solutions: $Z(x,y,s+2\pi)=Z(x,y,s)$ for all $s$. Introduce the functionals
\begin{equation}\label{BA-def}
\cB(Z) = \oint \fr \langle \bJ Z_s,Z\rangle \,\rd s \qand
\cA(Z) = \oint \fr \langle \bK Z_s,Z\rangle \,\rd s\,.
\end{equation}
Then a straightforward calculation confirms that
\begin{equation}\label{cwa}
\frac{\partial\ }{\partial x} \cB(Z) +
\frac{\partial\ }{\partial y} \cA(Z) = - \langle\nabla S(Z),\bM Z_t\rangle\,.
\end{equation}
Hence, loops of steady solutions have an exact conservation law with
components $(\cB,\cA)$.
It is the functions $(\cB,\cA)$, evaluated on the family of basic states,
that appear in the coefficients it the linearized
CN equation (\ref{steady-cn-linear}), in (\ref{kappa-def}),
and in the regularized characteristic CN equation (\ref{nrcn-eqn}).
However, to establish the connection between (\ref{cwa}) and (\ref{nrcn-eqn})
it is Noether's theorem applied to (\ref{MJKS-eqn}) that is key.  In
particular we have
\begin{equation}\label{noether-theory}
\bJ Z_s = \nabla \cB(Z) \qand \bK Z_s = \nabla \cA(Z)\,,
\end{equation}
where here the gradient is with respect to the inner product
\begin{equation}\label{ip-def}
  \langang Z,W\rangang =\oint\langle Z,W\rangle\,\rd s
  := \frac{1}{2\pi}\int_0^{2\pi}
  \langle Z,W\rangle\,\rd s\,.
\end{equation}
Later we will apply this theory to the basic periodic states,
with $s$ replaced by $\theta$ and averaging over $\theta$.  It is
the key identities (\ref{noether-theory}), with $Z_s$ replaced by
$\Zh_\theta$, that will connect the
values of the functionals $(\cB(k,\ell),\cA(k,\ell))$, and their
derivatives, to the structure matrices
$\bJ,\bK$ in the governing equations (\ref{MJKS-eqn}).

\subsection{The basic state}
\label{subsec-basicstate}

The basic state is a diagonal periodic pattern (roll) of the form
\begin{equation}\label{basicstate}
Z(x,y,t) = \widehat Z(\theta,k,\ell)\,,\quad \theta=kx+\ell y + \theta_0\,,
\end{equation}
with angle $\vartheta$ defined in (\ref{vartheta-def}),
and it satisfies
\begin{equation}\label{eq11}
  (k \bJ + \ell \bK)\Zh_\theta = \nabla S(\Zh)\,,\quad
  \mbox{for all}\ (k,\ell)\in\mathcal{D}\subset \R^2\,.
\end{equation}
This equation is a PDE if there is a cross section and
an ODE otherwise.  Solutions of (\ref{eq11}) are assumed to exist
and be smooth functions of $\theta,k,\ell$ for
all $(k,\ell)\in\mathcal{D}=\mathcal{D}^-\cup\mathcal{D}^+$.

Evaluating the functionals $\cB,\cA$ on this state and differentiating gives
\begin{equation}\label{eq10}
\begin{array}{rcl}
  && \cB_k = \langang\textbf{J}\widehat{Z}_\theta, \widehat{Z}_k \rangang\,,
  \qquad
  \cB_\ell = \langang\textbf{J}\widehat{Z}_\theta, \widehat{Z}_l \rangang\\[2mm]
  &&\cA_k = \langang\textbf{K}\widehat{Z}_\theta, \widehat{Z}_k \rangang\,,
  \qquad
  \cA_\ell = \langang\textbf{K}\widehat{Z}_\theta, \widehat{Z}_l \rangang\,.
  \end{array}
\end{equation}
The gradient structure assures that $\cB_\ell = \cA_k$.

\subsection{Linearization and solvability}
\label{subsec-linear}

Linearization of the governing equation (\ref{MJKS-eqn})
gives the following linear operator
\begin{align}
  \textbf{L} := D^2S(\Zh) - k\textbf{J} \frac{d\ }{d\theta} - \ell\textbf{K}\frac{d\ }{d\theta}\,.
\end{align}
We will use the same symbol ${\bf L}$ below when sketching
the derivation of the CN equation where the coefficients $(k,\ell)$
are replaced by $(k+q,\ell+r)$ respectively.
Differentiating \eqref{eq11} with respect to $\theta, k , \ell$ leads to equations that will be used in the modulation theory,
\begin{align}
    \textbf{L}\Zh_\theta = 0\,,\quad \textbf{L}\Zh_k = \textbf{J}\Zh_\theta\,,\quad \textbf{L}\Zh_\ell = \textbf{K}\Zh_\theta\,. \label{eq13}
\end{align}
The first of these shows that $\Zh_\theta$ is in the kernel of $\textbf{L}$, and we assume the kernel is no larger, so
\begin{equation}\label{L-kernel}
    \text{kernel}\{\textbf{L}\} = \text{span}\{\Zh_\theta\}.
\end{equation}
The other equations in \eqref{eq13} start the formation of two Jordan chains, one in the $\textbf{J}$ direction and one in the $\textbf{K}$ direction.
These Jordan chains will appear in the derivation of the CCN equation, and their
theory is developed in \S\ref{subsec-jordanchain}.

Throughout the modulation theory, and the construction of these Jordan chains,
inhomogeneous equations involving the linear operator $\textbf{L}$ will appear so a solvability condition is required. That is, for equations of the form $\textbf{L}V = F $ we have
\begin{align}
    \textbf{L}V = F \ \ \text{is solvable if and only if} \ \ \langang \hat{Z}_\theta, F\rangang = 0. \label{eq15}
\end{align}

\section{Re-appraisal of the coefficients in the CN equation}
\label{sec-cn-eqn-sketch}
\setcounter{equation}{0}

In this section a sketch of the derivation of the CN equation is given.
The sole aim to show how the coefficents in (\ref{steady-cn-linear}) can
be expressed in terms of derivatives of $\cB$ and $\cA$.

Given the basic state $\Zh$ in (\ref{basicstate}) the CN ansatz is
\begin{equation}\label{cn-ansatz-1}
    Z(x,y,t) = \Zh(\theta + \eps^{-1}\phi, k + q, \ell + r) + \eps W(\theta + \eps^{-1}\phi, X, Y,T)\,,
\end{equation}
    with $X=\eps x$, $Y=\eps y$ and $T=\eps^2t$, and $\phi,q,r$ all functions
    of $X,Y,T,\eps$.  With these scalings,
    \begin{equation}\label{eq18}
    \phi_X = q\,,\quad \phi_Y = r\,,\qand r_X = q_Y\,.
    \end{equation}
    Expand everything in a Taylor series, e.g.\ $W=W_1+\mathcal{O}(\eps)$,
    substitute into (\ref{MJKS-eqn}) and solve
    order by order in $\eps$.  At zeroth order the governing equation for
    the basic pattern (with $k\mapsto k+q$ and $\ell\mapsto \ell+r$)
    is recovered, and at first order we obtain
   \[
   \textbf{L}W_1 = \phi_T\bM\Zh_\theta+\bJ(q_X \Zh_k + r_X \Zh_\ell) + \textbf{K}(q_Y\Zh_{k} + r_Y \Zh_\ell)\,.
   \]
   Applying the solvability condition \eqref{eq15} then gives
\begin{align*}
  \phi_T\langang\Zh_\theta,\bM\Zh_\theta\rangang +
  q_X\langang\hat{Z}_\theta, \textbf{J}\hat{Z}_k \rangang + r_X \langang\hat{Z}_\theta, \textbf{J}\hat{Z}_l \rangang +  q_Y\langang\hat{Z}_\theta, \textbf{K}\hat{Z}_k \rangang + r_Y\langang \hat{Z}_\theta, \textbf{K}\hat{Z}_l \rangang  = 0\,.
\end{align*}
Taking into account that $\Zh$ is a function of $k+q$ and $\ell+r$, this
equation is in fact nonlinear and of the form
\[
\tau\phi_T = \frac{\partial\ }{\partial X}\cB(k+q,\ell+r) +
\frac{\partial\ }{\partial Y}\cA(k+q,\ell+r)\,,
\]
with $\tau =\langang\Zh_\theta,\bM\Zh_\theta\rangang$.
Linearizing and replacing $q=\phi_X$ and $r=\phi_Y$ then gives the
linearized CN or phase diffusion equation in standard form
\begin{equation}\label{steady-cn-linear-1}
  \tau\phi_T =
  \cB_k \phi_{XX} + (\cB_\ell+\cA_k)\phi_{XY} +
  \cA_{\ell}\phi_{YY}  \,,
\end{equation}
confirming \eqref{steady-cn-linear}.

\section{Derivation of the characteristic Cross-Newell equation}
\label{rccn-eqn-derivation}
\setcounter{equation}{0}

To derive the characteristic CN equation (\ref{nrcn-eqn}),
we start with the scaled independent variables
\begin{equation}\label{XYT-def-1}
X = \eps(x + C^{\pm} y)\,,\quad Y=\eps^3 y\,,\quad\mbox{and}\quad
T = \eps^4t\,.
\end{equation}
The revised modulation ansatz is
\begin{equation}\label{eq5.1}
    Z(x,y,t)  = \Zh(\theta + \eps \phi, k + \eps^2 q, \ell + \eps^2C^{\pm}q + \eps^4r) + \eps^3W(\theta, X, Y, T)\,,
\end{equation}
with $\phi,q,r$ all functions of $X,Y,T$, and the relations
\begin{equation}\label{eq4.3}
      \phi_X = q, \ \phi_Y = r, \ \phi_{XY} = q_Y = r_X\,.
\end{equation}
The novelty here, in addition to alignment with the characteristic
direction, is the mixed modulation of $k$ and $\ell$ in (\ref{eq5.1})
and this is inspired by the characteristic modulation theory in
\cite{r19}. It follows from (\ref{eq5.1}) that the precise form
of the matrix ${\bf C}(\eps)$ in (\ref{ansatz-2}) is
\begin{equation}\label{C-eps-def}
  {\bf C}(\eps) = \left[ \begin{matrix} 1 & 0 \\ C^{\pm} & \eps^2\end{matrix}
      \right]\,.
\end{equation}
In (\ref{eq5.1}), the
  phase perturbation $\theta\mapsto\theta+\eps\phi$ has been slowed down,
  and so everything can be expanded in a regular Taylor series.  The
  $W$ expansion is
\[
\eps^3W(\theta,X,Y,T) = \eps^3 W_3 +\eps^4 W_4 + \eps^5 W_5 +\cdots\,.
\]
Substitute into (\ref{MJKS-eqn}) and consider the $\eps^n$ terms for
each order $n$.  The $\eps^0$ equation recovers the
equation for $\Zh$, the first order equation recovers the kernel
condition (\ref{L-kernel}), and the second order equation recovers
the identity $q=\phi_X$.  The equations start to get interesting at
third order.

After some simplification, the terms remaining at this order are
\begin{equation}\label{W3-eqn}
    \textbf{L}W_3 = q_X (\textbf{J}+C\textbf{K})(\Zh_k + C\Zh_\ell).
\end{equation}
Applying solvability gives
\begin{align}
    \cB_k + C\cB_l + C \cA_k + C^2 \cA_l = 0\,,
\end{align}
showing that this equation is solvable if and only if $C=C^{\pm}$, that is,
one of the characteristics (\ref{C-plusminus-def}).

It will be shown below, in \S\ref{subsec-jordanchain},
that solvability of this equation generates a
Jordan chain of length four associated with (\ref{eq13}),
with the third generalized eigenvector denoted by $\xi_3$.
Thus, the solution for $W_3$ is
\begin{align}
    W_3 = q_X\xi_3\quad\mbox{(mod Ker({\bf L}))}\,.
\end{align}

The terms at fourth order in $\eps$ can be condensed into the form
\begin{equation}\label{W4-eqn}
    \textbf{L}\big(W_4 - \phi q_X(\xi_3)_\theta\big) = \phi_Y\textbf{K}\Zh_\theta - r\textbf{L}\Zh_\ell + q_{XX}(\textbf{J}+C^{\pm}\textbf{K})\xi_3\,.
\end{equation}
The first two terms on the right hand side cancel by using \eqref{eq4.3} and \eqref{eq13}. Now the equation is of the form of the
twisted Jordan chain (see \S\ref{subsec-jordanchain}) and
so the solution is generated by the fourth generalized eigenvector,
$\xi_4$, giving
\begin{align}
    W_4 = \phi q_X (\xi_3)_\theta + q_{XX}\xi_4\quad\mbox{(mod Ker({\bf L}))}\,.
\end{align}

\subsection{Fifth order}
\label{subsec-fifthorder}

The equation at fifth order in $\eps$ is where the solvability condition
produces the new equation (\ref{nrcn-eqn}).  Condensing terms, at
fifth order we find
\begin{equation}\label{LW5}
\begin{array}{rcl}
  \textbf{L}\widetilde{W}_5 &=&
  q_Y(\textbf{K}\Zh_k + C^{\pm}\textbf{K}\Zh_\ell) + r_X(\textbf{J}\Zh_\ell + C^{\pm}\textbf{K}\Zh_l) \\[2mm]
  &&\hspace{1.0cm} + qq_X((\textbf{J}+C^{\pm}\textbf{K})(\Zh_{kk} +2C^{\pm}\Zh_{\ell k} + (C^{\pm})^2\Zh_{\ell\ell} + (\xi_3)_\theta) - D^3S(\Zh)(\Zh_k + C^{\pm}\Zh_\ell, \xi_3))\\[2mm]
&&\hspace{2.0cm}    + q_{XXX}(\textbf{J}+C^{\pm}\textbf{K})\xi_4 + \phi_T\bM\Zh_\theta\,,
\end{array}
\end{equation}
where $\widetilde{W}_5$ contains $W_5$ and all other terms that can be shown to be in the range of $\textbf{L}$. The explicit form of $\widetilde{W}_5$ is not needed as it will vanish identically when the solvability condition is applied.

Using \eqref{eq4.3}, $r_X=q_Y$, and so the first two terms combine.
Then application of solvability to this term gives
\begin{align*}
    q_Y \ \text{coefficient:} \ \ \langang \Zh_\theta, \textbf{J}\Zh_l + \textbf{K}\Zh_k + 2C^{\pm}\textbf{K}\Zh_l \rangang = \cA_k + \cB_\ell + 2C^{\pm}\cA_\ell\,.
\end{align*}
Using the formula for characteristics (\ref{C-plusminus-def}), 
\begin{align*}
   \cA_k + \cB_\ell + 2C^{\pm}\cA_\ell = \mp 2 \sqrt{-\Delta_{zz}}\,.
\end{align*}
The solvability condition for the $qq_X$ term is
\begin{align*}
 -\kappa := \langang \Zh_\theta, (\textbf{J}+C^{\pm}\textbf{K})(\Zh_{kk} +2C^{\pm}\Zh_{lk} + (C^{\pm})^2\Zh_{ll}
  + (\xi_3)_\theta)- D^3S(\Zh)(\Zh_k + C^{\pm}\Zh_l, \xi_3) \rangang \,.
\end{align*}
Calculating, and using the theory in \cite{r19}, gives
\[
\kappa = (\partial_k + C^{\pm}\partial_\ell)^2(\cB + C^{\pm}\cA)\,,
\]
where the derivatives are taken with $C^{\pm}$ fixed.
The coefficient of regularization is
\begin{align*}
    q_{XXX} \ \text{coefficient:} \ \    \langang \Zh_\theta, (\textbf{J}+C^{\pm}\textbf{K})\xi_4 \rangang = - \mathscr{K}.
\end{align*}
At this point the coefficient is just called $-\mathscr{K}$, but we will see
in the next section that this coefficient emerges from the twisted
Jordan chain theory.  The coefficient of the $\phi_T$ term is the
same $\tau$ term as in the classical CN theory, the only difference being
the definition of $T$.

Putting this all together, and re-introducing $\phi$ using
$q=\phi_X$, we obtain
\begin{align}
\tau \phi_T =    \pm 2 \sqrt{-\Delta_{zz}}\,\phi_{XY} + \kappa \phi_X\phi_{XX} + \mathscr{K}\phi_{XXXX} \,.
\end{align}
This completes the derivation of (\ref{nrcn-eqn}).

\subsection{A twisted Jordan chain}
\label{subsec-jordanchain}

There are two Jordan chains of length two in the linearization about a
periodic pattern.  They were introduced in \eqref{eq13}.
However it is the combination of these two chains, what we call a
twisted Jordan chain, that appears in the
derivation of (\ref{nrcn-eqn}). The first step in the chain is
still ${\bf L}\Zh_\theta=0$, but the second step is
\[
  {\bf L}(\Zh_k + C^{\pm}\Zh_\ell) = (\bJ + C^{\pm}\bK)\Zh_\theta\,.
  \]
  We will show that this Jordan chain continues to length four.
  To formalize the chain, introduce
  the notation $\{\xi_1,\xi_2,\xi_3,\xi_4\}$ with
  \[
  \xi_1 = \Zh_\theta\,,\quad \xi_2 =\Zh_k+C^{\pm}\Zh_\ell\quad\mbox{and}\quad
  {\bf L}\xi_2 = (\bJ+C^{\pm}\bK)\xi_1\,.
  \]
  The next element in the chain, if it exists, satisfies
  \begin{equation}\label{xi3-eqn}
    {\bf L}\xi_3 = (\bJ + C^{\pm}\bK)\xi_2\,.
  \end{equation}
    This equation is indeed solvable since
    \[
    \begin{array}{rcl}
    \langang \Zh_\theta,(\bJ + C^{\pm}\bK)\xi_2\rangang
    &=& \langang\Zh_\theta,(\bJ + C^{\pm}\bK)\xi_2\rangang\\[2mm]
    &=& \langang\Zh_\theta,(\bJ + C^{\pm}\bK)(\Zh_k+C^{\pm}\Zh_\ell)\rangang\\[2mm]
    &=& \langang\Zh_\theta,\bJ(\Zh_k+C^{\pm}\Zh_\ell)\rangang+
    C^{\pm}\langang\Zh_\theta,\bK(\Zh_k+C^{\pm}\Zh_\ell)\rangang\\[2mm]
    &=& -\cB_k - C^{\pm}\cB_\ell - C^{\pm}\cA_k - (C^{\pm})^2\cA_\ell\,.
    \end{array}
    \]
    This right-hand side vanishes identically when $C^{\pm}$ are
    characteristics (\ref{C-plusminus-def}).
    Hence the Jordan chain has length three.  It is immediate however that the
    length is four: there exists $\xi_4$ satisfying
  \begin{equation}\label{xi4-eqn}
      {\bf L}\xi_4 = (\bJ + C^{\pm}\bK)\xi_3\,.
  \end{equation}
      The existence of $\xi_4$ can be confirmed by calculating solvability or by
      noting that zero is an eigenvalue of multiplicity four of
      the operator pencil ${\bf L}-\lambda(\bJ+C^{\pm}\bK)$, since ${\bf L}$
      is symmetric and $\bJ+C^{\pm}\bK$ is skew-symmetric.

      The $\xi_3$ chain equation (\ref{xi3-eqn})
      is of central importance in the third-order
      equation (\ref{W3-eqn}), and the $\xi_4$ chain (\ref{xi4-eqn})
      if of central importance in the fourth-order equation (\ref{W4-eqn}).

      For the chain to terminate at four, the equation
      \[
        {\bf L}\xi_5 = (\bJ+C^{\pm}\bK)\xi_4\,,
        \]
        should be not solvable; that is
        \[
        \langang \Zh_\theta, (\textbf{J}+C^{\pm}\textbf{K})\xi_4\rangang \neq 0\,.
          \]
          Call this right-hand side $-\mathscr{K}$.  It is this $\mathscr{K}$
          that appears as the ``regularizing'' coefficient of the fourth
          derivative term in (\ref{nrcn-eqn}).

\subsection{Localized solutions along the third direction}
\label{sec-solitons}

The steady part of (\ref{nrcn-eqn}) is the KdV equation in the wavenumber
$q=\phi_X$,
\begin{align}
    q_Y + \frac{(\partial_k + c\partial_l)^2(\cB + c \cA)}{\pm 2 \sqrt{\Delta_L}}qq_X + \frac{\mathscr{K}}{\pm 2\sqrt{\Delta_L}} q_{XXX} = 0.
\end{align}
Rescale variables to put the KdV in standard form; let
\begin{align}
    X = \alpha \widetilde{X}, \ Y = \beta \widetilde{Y}, \ q = \gamma \widetilde{q}\,,
\end{align}
and choose $\alpha, \beta, \gamma$ such that
\begin{align}
     \frac{\gamma}{\beta} = \frac{\gamma^2}{\alpha}\frac{(\partial_k + c\partial_l)^2(\cB + c \cA)}{\pm 2 \sqrt{\Delta_L}} = \frac{\gamma}{\alpha^3}\frac{\mathscr{K}}{\pm 2\sqrt{\Delta_L}}\,,
\end{align}
then we have
\begin{align}
    \widetilde{q}_{\widetilde{Y}} + \widetilde{q}\widetilde{q}_{\widetilde{X}} + \widetilde{q}_{\widetilde{X}\widetilde{X}\widetilde{X}} = 0.
\end{align}
This equation has the standard KdV solitary wave.  However, it is not dynamic,
it just appears at some angle.  This is the third angle of the theory.
Let
\[
\zeta = \widetilde X - c_3 \widetilde Y\qand \widetilde q(\widetilde X,\widetilde Y) = \widehat q(\zeta)\,.
\]
Then $\widehat q$ satisfies the usual ODE
$- c_3 \widehat{q}' + \widehat{q}\widehat{q}' + \widehat{q}''' = 0$
with solitary wave
\begin{align}
    \widetilde q(\widetilde X,\widetilde Y) = 3c_3 \,{\rm sech}^2\left(\frac{\sqrt{c_3}}{2}(\widetilde X-c_3\widetilde Y)\right)\,.
\end{align}
This ``solitary wave'' is bi-asymptotic to rolls as $\zeta\to\pm\infty$.

\section{Example: Ginzburg-Landau equation}
\label{sec-gl}
\setcounter{equation}{0}

Consider the real Ginzburg-Landau (RGL) equation in two space
dimensions, 
\begin{equation}\label{rgl-eqn}
  \Psi_t = \Psi_{xx}+\Psi_{yy} + \Psi -|\Psi|^2\Psi \,,
\end{equation}
for the complex-valued function $\Psi(x,y,t)$. There are
a number of ways to write the equation in the multisymplectic
form (\ref{MJKS-eqn}).  Let
\[
\bu=\begin{pmatrix} u_1\\ u_2\end{pmatrix}\quad\mbox{with}\quad
u_1+\ri u_2 = \Psi\,,
\]
and introduce derivatives $\bp=\bu_x$ and $\br=\bu_y$.  Then
the RGL equation in real coordinates is
\[
\underbrace{\left[\begin{matrix} {\bf I} & {\bf 0} & {\bf 0} \\
    {\bf 0} & {\bf 0} & {\bf 0} \\
    {\bf 0} & {\bf 0} & {\bf 0} \end{matrix}\right]}_{\bM}
\begin{pmatrix}    \bu\\ \bp\\ \br\end{pmatrix}_t
+
  \underbrace{\left[\begin{matrix} {\bf 0} & -{\bf I} & {\bf 0} \\
    {\bf I} & {\bf 0} & {\bf 0} \\
    {\bf 0} & {\bf 0} & {\bf 0} \end{matrix}\right]}_{\bJ}
  \begin{pmatrix}    \bu\\ \bp\\ \br\end{pmatrix}_x+
    \underbrace{\left[\begin{matrix} {\bf 0} & {\bf 0} & -{\bf I} \\
    {\bf 0} & {\bf 0} & {\bf 0} \\
    {\bf I} & {\bf 0} & {\bf 0} \end{matrix}\right]}_{\bK}
    \begin{pmatrix}    \bu\\ \bp\\ \br\end{pmatrix}_y =
      \begin{pmatrix} \delta S/\delta\bu\\ \delta S/\delta\bp\\
        \delta S/\delta \br \end{pmatrix}\,,
      \]
      with
      \[
      S = \fr \|\textbf{u}\|^2 + \fr \|\textbf{p}\|^2 + \fr \|\textbf{r}\|^2 - \fr \|\textbf{u}\|^4\,.
      \]
      However a more elegant formulation is obtained by adding the
      cross derivative equation $\bp_y=\br_x$.  This requires a new
      coordinate $\bw$.  $S(Z)$ remains the same, although we now think
      of $Z$ as having real dimension $8$.  The matrix $\bM$ is the same
      with an extra vector row and column of zeros, and the matrices
      $\bJ$ and $\bK$ are now
\begin{align}
    &Z = \begin{pmatrix}\textbf{u} \\ \textbf{p} \\ \textbf{r} \\ \textbf{w} \end{pmatrix}, \ \textbf{J} = \begin{pmatrix}\textbf{0} & -\textbf{I} & \textbf{0} & \textbf{0} \\ \textbf{I} & \textbf{0} & \textbf{0} & \textbf{0} \\ \textbf{0} & \textbf{0} & \textbf{0} & \textbf{I} \\ \textbf{0} & \textbf{0} & -\textbf{I} & \textbf{0}\end{pmatrix} , \ \textbf{K} = \begin{pmatrix}\textbf{0} & \textbf{0} & -\textbf{I} & \textbf{0} \\ \textbf{0} & \textbf{0} & \textbf{0} & -\textbf{I} \\ \textbf{I} & \textbf{0} & \textbf{0} & \textbf{0} \\ \textbf{0} & \textbf{I} & \textbf{0} & \textbf{0} \end{pmatrix}\,. 
\end{align}
The structure matrices $\bJ$ and $\bK$ are now invertible, and each generate
standard symplectic structures.  Moreover, the three matrices $\{\bJ,\bK,\bJ\bK\}$ generate a representation of
the quaternions, although that property will not be needed here.

The basic state is the family of rolls $\Psi=\widehat\Psi\re^{\ri\theta}$,
with $\theta=kx+\ell y$, which in real coordinates is
\[
\bu(x,y,t) = \Rg_\theta\widehat\bu(k,\ell)
\quad\mbox{with}\quad
\Rg_\theta:=\left[\begin{matrix} \cos\theta & -\sin\theta \\
    \sin\theta & \cos\theta \end{matrix}\right]\,,
\]
and the full $8-$dimensional representation
\begin{equation}\label{Zh-formula}
\Zh(\theta,k,\ell) =
   {G}_\theta \Zh(k,\ell)\,,\quad {G}_\theta = \Rg_\theta\oplus
\Rg_\theta\oplus\Rg_\theta\oplus\Rg_\theta\,,\quad
\Zh(k,\ell) = \begin{pmatrix} \widehat\bu(k,\ell)\\ k\bJ_2\widehat\bu(k,\ell)\\
  \ell\bJ_2\widehat\bu(k,\ell)\\ {\bf 0}\end{pmatrix}\,,
\end{equation}
where $\bJ_2$ is the action of $\bJ$ on $\R^2$.  Substituting this
representation
into the governing equations gives
\begin{equation}\label{klu-eqn}
k^2 + \ell^2 + \|\widehat\bu\|^2 = 1 \,.
\end{equation}
From this we get that the set of existence is the unit disc
\[
\mathcal{D}(k,\ell) =\big\{ (k,\ell)\in\R^2\ :\ k^2+\ell^2<1\big\}\,.
\]

\begin{figure}
    \centering
    \includegraphics[width=7cm]{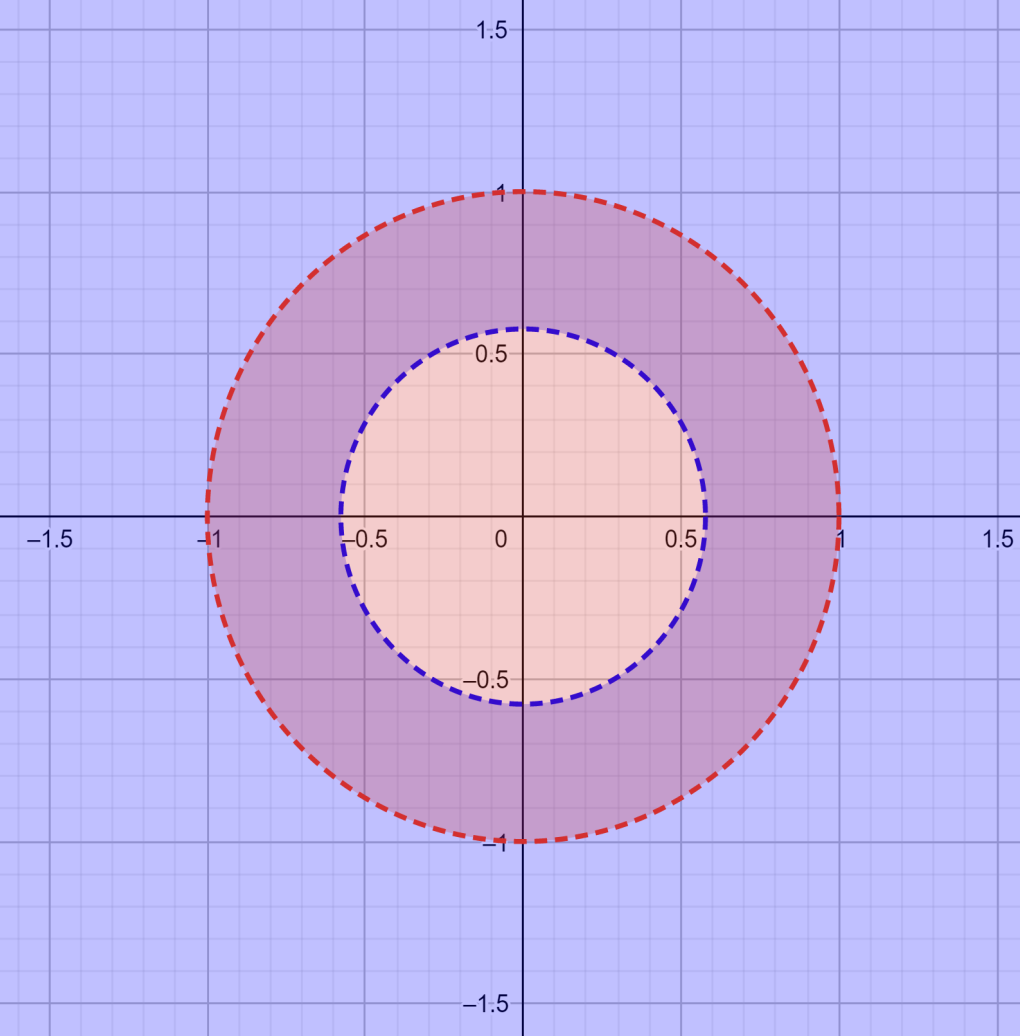}
    \caption{Existence and zig-zag regions in the $k,\ell$ plane.}
    \label{fig1}
\end{figure}

Evaluating the basic state on the components of the conservation law gives
\begin{equation}\label{rgl-AB-calc}
\cB(k,\ell) = k\|\widehat\bu\|^2 = k(1-k^2-\ell^2)\qand
   \cA = \ell\|\widehat\bu\|^2 = \ell(1-k^2-\ell^2)\,.
\end{equation}
The characteristics are
\[
C^{\pm} = -\frac{(\cB_\ell+\cA_k)}{2\cA_{\ell}}
\pm \frac{1}{\cA_{\ell}}\sqrt{-\Delta_{zz}}\,.
\]
From (\ref{rgl-AB-calc}) the terms are
\[
\cA_\ell = 1-k^2 - 3\ell^2\,,\quad
\cB_\ell=\cA_k = -2k\ell\,,
\]
and
\[
\Delta_{zz}={\rm det}\left[ \begin{matrix} \cB_k & \cB_\ell \\
    \cA_k & \cA_\ell \end{matrix}\right] =
{\rm det}\left[ \begin{matrix} 1-3k^2-\ell^2 & -2k\ell \\
    -2k\ell & 1-3\ell^2-k^2 \end{matrix}\right] =
(1-3k^2-3\ell^2)\|\widehat\bu\|^2\,.
\]
Real characteristics exist in the annulus
\[
\mathcal{D}^- =\big\{ (k,\ell)\in\mathcal{D}\ :\ \Delta_{zz}<0\big\} =
\big\{(k,\ell)\in\mathcal{D}\ :\ 1> k^2+\ell^2>\frac{1}{3}\big\}\,.
\]
The existence region (boundary in red) and zig-zag unstable region (inner boundary
in blue) are shown in Figure \ref{fig1}. 
The CCN equation is operational at all points in the annulus.  Since
$\cB$ and $\cA$ are explicitly known, it is straightfoward to compute
all the coefficients in the CCN equation.  The first two are
\begin{align}
    &\phi_{XY}\,: \ \ -\cA_k - \cB_\ell - 2C^{\pm}\cA_\ell = 4k\ell - 2C^{\pm}(1-3\ell^2-k^2) \\
    &\phi_X\phi_{XX}\,: \ \ (\partial_k + C^{\pm}\partial_\ell)^2(\cB +C^{\pm}\cA) = -6k - 6C^{\pm}\ell - (2k+4\ell)(C^{\pm})^2 - 6(C^{\pm})^3\ell\,.
\end{align}
To obtain an explicit expression for the coefficient of $\phi_{XXXX}$
requires computation of the twisted Jordan chain. This computation is
a bit lengthy but worth doing as it will tell us whether the
fourth derivative term is regularizing or destabilizing.

For the computation of $\mathscr{K}$, we need the operator $\textbf{L}$. It is
\begin{align}
    \textbf{L} = \begin{pmatrix} \textbf{A} & k\textbf{J}_2 & \ell\textbf{J}_2 & \textbf{0}_2 \\ -k\textbf{J}_2 & \textbf{I}_2 & \textbf{0}_2 & \ell\textbf{J}_2 \\ -\ell\textbf{J}_2 & \textbf{0} & \textbf{I}_2 & -k\textbf{J}_2 \\ \textbf{0}_2 &-\ell\textbf{J}_2 & k\textbf{J}_2 & \textbf{0}_2\end{pmatrix},
\end{align}
where
\begin{align}
  \textbf{A} = (k^2+\ell^2){\bf I}_2 - 2\widehat\bu\widehat\bu^T
  \,,\quad \textbf{J}_2 = \begin{pmatrix}0 & -1 \\1 & 0 \end{pmatrix}.
\end{align}
The first element in the Jordan chain is $\Zh_\theta$ which is obtained
by differentiating $\Zh$ in (\ref{Zh-formula}),
\begin{align}
    \Zh_\theta = G_\theta \begin{pmatrix}\textbf{J}_2 \widehat\bu \\ -k \widehat\bu \\ -l \widehat\bu \\ \textbf{0} \end{pmatrix}.
\end{align}
To obtain $\xi_2$ we need $\Zh_k$ and $\Zh_\ell$.  They are
\begin{align}
  \Zh_k  = G_\theta \begin{pmatrix}-\frac{k}{\|\widehat\bu\|^2}\widehat\bu \\ \frac{\|\widehat\bu\|^2 - k^2}{\|\widehat\bu\|^2}\textbf{J}_2 \widehat\bu \\-\frac{k\ell}{\|\widehat\bu\|^2}\textbf{J}_2 \widehat\bu \\ \textbf{0} \end{pmatrix}
  \qand
\Zh_\ell = G_\theta \begin{pmatrix}-\frac{l}{\|\widehat\bu\|^2} \widehat\bu \\ -\frac{k\ell}{\|\widehat\bu\|^2}\textbf{J}_2 \widehat\bu \\ \frac{\|\widehat\bu\|^2 - \ell^2}{\|\widehat\bu\|^2}\textbf{J}_2 \widehat\bu \\ \textbf{0}  \end{pmatrix}\,.
\end{align}
To obtain $\xi_3$ we need to
solve $\textbf{L}\xi_3 = (\textbf{J} + C^{\pm}\textbf{K})(\Zh_k + C^{\pm}\Zh_\ell)$, 
\begin{align}
    \xi_3 = G_\theta\begin{pmatrix}\textbf{0} \\ \frac{-k-C^{\pm}\ell}{\|\widehat\bu\|^2}\widehat\bu \\ \frac{-C^{\pm}k - (C^{\pm})^2\ell}{\|\widehat\bu\|^2}\widehat\bu \\ \textbf{0} \end{pmatrix}.
\end{align}
The expressions for $\xi_2$ and $\xi_3$ are particular solutions, and
an arbitrary amount of homogeneous solution (kernel of ${\bf L}$) can be
added.  However, the homogeneous terms cancel out in the calculation of
$\mathscr{K}$.

We do not need to calculate $\xi_4$.  Using the Jordan chain
(\ref{subsec-jordanchain}), we have
\begin{align*}
    \mathscr{K} &=  -\langang\Zh_\theta, (\textbf{J}+C^{\pm}\textbf{K})\xi_4 \rangang \\
    &=\langang(\textbf{J}+C^{\pm}\textbf{K}) \xi_1, \xi_4 \rangang \\
    &= \langang \textbf{L}\xi_2, \xi_4 \rangang \\
    &= \langang \xi_2, \textbf{L}\xi_4 \rangang \\
    &= \langang \xi_2, (\textbf{J} + C^{\pm} \textbf{K})\xi_3 \rangang.
\end{align*}
Substituting $\xi_2=\Zh_k+C^{\pm}\Zh_\ell$ and $\xi_3$ into this formula
then gives
\begin{align}
    \mathscr{K} = - \frac{1}{\|\widehat\bu\|^2}(1+(C^{\pm})^2)(k+C^{\pm}l)^2\,,
\end{align}
which is strictly negative and therefore stabilizing.  The coefficient
of the $\phi_T$ term is
\[
\tau = \langang \bM\Zh_\theta,\Zh_\theta\rangang = \|\widehat\bu\|^2\,.
\]
Hence, the CCN equation for modulation of the rolls along characteristic
directions is then
\begin{equation}\label{RGL-CCN}
\begin{array}{rcl}
  &&(1-k^2-\ell^2)\phi_T = \big(4k\ell - 2C^{\pm}(1-3\ell^2-k^2)\big)\phi_{XY} \\[2mm]
  &&\hspace{3.0cm}
  + \Big(-6k - 6C^{\pm}\ell - (2k+4\ell)(C^{\pm})^2 - 6(C^{\pm})^3\ell\Big)\phi_X\phi_{XX} \\[2mm]
  &&\hspace{6.0cm} - \frac{1}{\|\widehat\bu\|^2}\Big(1+(C^{\pm})^2\Big)(k+C^{\pm}\ell)^2\phi_{XXXX} = 0\,.
    \end{array}
\end{equation}
This equation is asymptotically valid in the annulus in Figure \ref{fig1}.
The outer boundary of the annulus is determined by the region of existence.
The inner boundary, the blue border in Figure \ref{fig1} is the curve along which
$\Delta_{zz}=0$.  This singularity corresponds to coalescing characteristics.
Along this curve there is only one characteristic and the theory breaks down.
However one can again re-modulate.  A theory for coalescing characteristics
in Whitham theory for time evolution equations has been developed in \cite{br17}
and could be adapted to this setting.  Based on this theory, it
is expected that re-modulation will generate
a two-way Boussinesq
equation, leading to a replacement of the $\phi_{XY}$ term
with a $\phi_{XYY}$ term.

\section{Concluding remarks}
\label{sec-cr}
\setcounter{equation}{0}

The theory is valid in an open set of wavenumbers, $\mathcal{D}^-$,
determined by the intersection of the set of existence $\mathcal{D}(k,\ell)$ with the
set of wavenumbers where $\Delta_{zz}<0$.  This set is a large open set with
the same Lebesgue measure as the existence set $\mathcal{D}$.
In the RGL example it is an annulus.
If we can cover the complement, $\mathcal{D}^+$, where $\Delta_{zz}>0$, then the
CCN equation would be asymptotically valid in the same region as the
CN equation.  Hence the most intriguing potential new direction is the case
where the characteristics are complex.
Mathematically it looks possible, $X=\eps(x+C^{\pm}y)$ would be complex
and the coefficients in the CCN equation would be complex, but the details
of the complexification of the Jordan chain, and the physical interpretation of
a complex CCN equation are not clear.

\hfill{\textsf{\today}}

\end{document}